\RequirePackage{fix-cm}
\documentclass[smallextended]{svjour3}       % onecolumn (second format)
\smartqed  % flush right qed marks, e.g. at end of proof
\usepackage{amssymb,amsmath,bm,hyperref}
\usepackage{amssymb,amsmath,bm,geometry,graphics,graphicx,url,color}
\def\pmatrix{\left(\begin{array}}
\def\endpmatrix{\end{array}\right)}
\def\RR{\mathbb{R}}

\def\I{{\cal I}}

\def\P{{\cal P}}

\def\dd{\mathrm{d}}

\newtheorem{theo}{Theorem}

\newtheorem{rem}{Remark}
\newtheorem{defi}{Definition}

\def\bfb{{\bm{b}}}
\def\bfc{{\bm{c}}}

\def\bfe{{\bm{e}}}

\def\bfp{{\bm{p}}}
\def\bfq{{\bm{q}}}

\def\bfw{{\bm{w}}}

\def\bfy{{\bm{y}}}

\def\bfgamma{{\bm{\gamma}}}

\begin{document}

\title{Spectrally accurate space-time solution of Hamiltonian PDEs}

\author{Luigi Brugnano  \and  Felice Iavernaro  \and Juan I. Montijano \and Luis R\'andez}

\institute{L. Brugnano \at
Dipartimento di Matematica e Informatica ``U.\,Dini'', Universit\`a di Firenze, Viale Morgagni 67/A, 50134 Firenze, Italy. \\
\email{luigi.brugnano@unifi.it}
\and
F. Iavernaro \at 
Dipartimento di Matematica, Universit\`a di Bari, Via Orabona 4, 70125 Bari, Italy. \\
\email{felice.iavernaro@uniba.it}
\and 
J.\,I.\,Montijano \at 
Departamento  de Matem\'{a}tica Aplicada, Universidad de Zaragoza, Pza. San Francisco s/n,
50009 Zaragoza, Spain. \\
\email{monti@unizar.es}
\and 
L. R\'andez \at 
Departamento  de Matem\'{a}tica Aplicada, Universidad de Zaragoza, Pza. San Francisco s/n,
50009 Zaragoza, Spain. \\
\email{randez@unizar.es}}

\dedication{Dedicated to John Butcher, on the occasion of his 85-th birthday}

\date{}

\maketitle

\begin{abstract} Recently, the numerical solution of multi-frequency, highly-oscillatory Hamiltonian problems has been attacked by using Hamiltonian Boundary Value Methods (HBVMs) as spectral methods in time.
When the problem derives from the space semi-discretization of (possibly Hamiltonian) partial differential equations (PDEs), the resulting problem may be {\em stiffly-oscillatory}, rather than highly-oscillatory. In such a case, a different implementation of the methods is needed, in order to gain the maximum efficiency.



\end{abstract}

\section{Introduction}\label{intro}

Multi-frequency highly-oscillatory problems have been recently attacked by using Hamiltonian Boundary Value Methods (HBVMs) as spectral methods in time \cite{BMR2018}. The proposed approach has proven to be very efficient when solving a number of severe highly-oscillatory problems, allowing to effectively and accurately ``resolve'' all high-frequency components in the solution. Sometimes, however, the problem is only {\em stiffly-oscillatory}, rather than highly-oscillatory, which means that the high-frequency components in the solution all have a {\em very small amplitude}. This is the case, for example, of problems deriving from the space semi-discretization of time-dependent PDEs having a relatively smooth solution. 

In this paper, we shall consider Hamiltonian PDEs with periodic boundary conditions possessing a soliton-type solution. In such a case, the two implementing criteria devised in \cite{BMR2018}, which are aimed at grasping {\em all} the high-frequencies, may be too much stringent and, therefore, the approach could become less efficient.
In fact, the highest frequencies with negligible amplitude, could be more conveniently omitted, since their contribution to the accuracy of the solution is actually marginal. In order to restore the efficiency of HBVMs used as spectral methods, we here propose an adaptive implementation of the methods, able to overcome this drawback, still providing a {\em practical} spectral accuracy in time. Coupling this approach with a spectrally accurate space semi-discretization will result in a {\em spectrally accurate space-time solution} of Hamiltonian PDEs.

With this premise, the structure of the paper is as follows: in Section~\ref{key} we describe the main differences with the approach described in \cite{BMR2018}; in Section~\ref{space} we provide some details on the semi-discrete problem derived from the space discretization of the considered Hamiltonian PDEs; in Section~\ref{numtest} some numerical tests are reported; at last, a few conclusions are given in Section~\ref{fine}.

\section{Basic facts}\label{key}
We are here concerned with the numerical solution of Hamiltonian problems in the form
\begin{equation}\label{odep}
\dot y = J\left[ A y + \nabla f(y) \right] \equiv \phi(y), \qquad y(0)=y_0\in\RR^{2m}, \qquad J = \pmatrix{cc} &I_m\\ -I_m\endpmatrix,
\end{equation}
where, in general, $I_r\in\RR^{r\times r}$ is the identity matrix, $A$ is a symmetric and positive semi-definite matrix such that, in a neighbourhood of  the solution,
\begin{equation}\label{hop}
\|A\| \equiv \omega \gg \| \nabla f\|,
\end{equation}
where $\|\cdot\|$ denotes the 2-norm, so that $\omega=\rho(A)$. 
%even though any norm could in principle be used, thus providing an upper bound for $\rho(A)$. 
Moreover, hereafter we assume $f$ to be suitably smooth, e.g., analytical. Problem (\ref{odep}) is clearly Hamiltonian with Hamiltonian
\begin{equation}\label{H}
H(y) = \frac{1}2 y^\top Ay + f(y) \qquad\qquad\left(~\Rightarrow \phi(y)=J\nabla H(y)~\right).
\end{equation}
A differential system in the form (\ref{odep})--(\ref{hop}) provides an instance of a, possibly multi-frequency,\footnote{Depending on the occurrence of different large eigenvalues of matrix $A$.} highly-oscillatory problems. We refer to the recent monograph \cite{WW2018} for an account of the various approaches used so far for dealing with such problems.

More recently, in \cite{BMR2018} a spectral method along the Legendre polynomial basis,
\begin{equation}\label{orto}
\deg(P_i)=i, \qquad \int_0^1 P_i(x)P_j(x)\dd x = \delta_{ij}, \qquad \forall i,j=0,1,\dots,
\end{equation}
has been defined, as is sketched below, based on the approach defined in \cite{BIT2012}. We start considering the expansion of the right-hand side of the differential equation in (\ref{odep}), on the interval $[0,h]$, along the orthonormal basis (\ref{orto}):
\begin{equation}\label{expfi}
\dot y(ch) = \sum_{j\ge0} P_j(c) \gamma_j, \qquad c\in[0,1], \qquad \gamma_j = \int_0^1 P_j(\tau)\phi(y(\tau h))\dd\tau, \qquad j=0,1,\dots.
\end{equation}
Integrating term by term the first equation in (\ref{expfi}), and imposing the initial condition in (\ref{odep}), then gives:
\begin{equation}\label{expy}
y(ch) = y_0 + h\sum_{j\ge0} \int_0^c P_j(x)\dd x\, \gamma_j, \qquad c\in[0,1].
\end{equation}
At this point, two facts have to be taken into account, in order for (\ref{expy}) to become an effective method, when using a finite precision arithmetic with machine epsilon $u$:\footnote{E.g., $u\approx 10^{-16}$, for the double precision IEEE.}
\begin{enumerate}
\item in the expansion (\ref{expy}) the coefficients $\gamma_j$ which are too small can be neglected. By considering that their norm is approximately decreasing from a certain index on, one has that
\begin{equation}\label{expys}
y(ch) \doteq y_0 + h\sum_{j=0}^{s-1} \int_0^c P_j(x)\dd x\, \gamma_j, \qquad c\in[0,1],
\end{equation}
where $\doteq$ stands for {\em equal within machine precision}, with the index $s$ satisfying
\begin{equation}\label{ratio}
\rho_s \equiv \frac{\|\gamma_{s-1}\|}{\max_{i=0,\dots,s-1} \|\gamma_i\|} ~\le~u\,;
\end{equation}

\item with reference to the coefficients $\gamma_j$ defined in (\ref{expfi}), one has
\begin{equation}\label{quadra}
\gamma_j \doteq \sum_{i=1}^k b_i P_j(c_i)\phi( y(c_ih)),
\end{equation}
by using an enough accurate quadrature which, hereafter, we choose as the Gauss-Legendre quadrature formula of order $2k$, for a suitable $k>s$. Unless better choices are available, we consider the choice used in \cite[Eq.\,(34)]{BMR2018}, i.e.,
\begin{equation}\label{kappa}
k = \max\{ s+2, 20\}.
\end{equation}
\end{enumerate}
Next, by considering that:
\begin{itemize}
\item by setting $Y_i \doteq y(c_ih)$, from (\ref{expys}) and (\ref{quadra}) one obtains
$$Y_i = y_0 + h\sum_{j=1}^k \left[ b_j \sum_{\ell=0}^{s-1} \int_0^{c_i} P_\ell(x)\dd x\, P_\ell(c_j)\right] \phi(Y_j), \qquad i=1,\dots,k,$$

\item the new approximation, by virtue of the orthogonality conditions (\ref{orto}), is defined as
\begin{equation}\label{yh}
y(h) \doteq y_1 = y_0 + h\sum_{i=1}^k b_i \phi( Y_i),
\end{equation}
\end{itemize}
one eventually arrives at the $k$-stage Runge-Kutta method defined by the following Butcher tableau,
\begin{equation}\label{butab}
\begin{array}{c|c}
\bfc & \I_s\P_s^\top\Omega \\ \hline \\[-3mm] &\bfb^\top
\end{array}\,,
\end{equation}
with
\begin{eqnarray}\nonumber
\bfc &=& (c_1,\dots,c_k)^\top, \quad \bfb~=~(b_1,\dots,b_k)^\top, \qquad \Omega = \pmatrix{ccc} b_1\\ &\ddots\\ &&b_k\endpmatrix,\\
\label{butab1}\\ \nonumber
\I_s &=& \pmatrix{ccc} \int_0^{c_1} P_0(x)\dd x& \dots &\int_0^{c_1} P_{s-1}(x)\dd x\\
\vdots & &\vdots\\
\int_0^{c_k} P_0(x)\dd x& \dots &\int_0^{c_k} P_{s-1}(x)\dd x\endpmatrix, \quad
\P_s ~=~ \pmatrix{ccc} P_0(c_1) &\dots &P_{s-1}(c_1)\\
\vdots & &\vdots\\ P_0(c_k) &\dots &P_{s-1}(c_k)\endpmatrix.
\end{eqnarray}

\medskip
\begin{defi}\label{hbvmks} The $k$-stage Runge-Kutta method (\ref{butab})--(\ref{butab1}) is called {\em Hamiltonian Boundary Value Method} with parameters $(k,s)$. In short, HBVM$(k,s)$.\end{defi}

It is quite clear that, by choosing $k$ large enough, the method is able to conserve, either exactly or within the round-off error level, the Hamiltonian $H$ along the numerical solution. We also mention that the HBVM$(k,s)$ family may be thought of as a generalization of Gauss-Legendre collocation methods in that HBVM$(s,s)$ coincides with the $s$-stage Gauss integrator.
We refer to the monograph \cite{LIMbook2016}, and to the recent review paper \cite{BI2018}, for full details about HBVMs.

Moreover, we observe that, when the parameters $k$ and $s$ are chosen such that (\ref{expys}) and (\ref{quadra}) hold true, they provide a spectrally accurate in time method for the solution of (\ref{odep})--(\ref{hop}).

It is also worth mentioning that the discrete problem generated by a HBVM$(k,s)$ method can be cast in terms of the $s$ coefficients $\gamma_0,\dots,\gamma_{s-1}$ in (\ref{expfi}),  thus leading to a nonlinear system having (block) dimension $s$ {\em independently} of $k$ \cite{BIT2011}. In fact, by setting (see (\ref{H}))
$$\bfgamma = \pmatrix{c} \gamma_0\\ \vdots \\ \gamma_{s-1}\endpmatrix, \quad Y = \pmatrix{c} Y_1\\ \vdots \\ Y_k\endpmatrix, \quad \phi(Y) = \pmatrix{c} \phi(Y_1)\\ \vdots \\ \phi(Y_k)\endpmatrix \equiv (I_k\otimes J)\nabla H(Y), \quad \bfe = \pmatrix{c}1\\ \vdots \\ 1\endpmatrix \in\RR^k,$$
from (\ref{expys})--(\ref{quadra}) one obtains
$$Y = \bfe\otimes y_0 + h\I_s\otimes I_{2m} \bfgamma, \qquad \bfgamma = \P_s^\top\Omega \otimes J\, \nabla H(Y),$$
which, combined together, provide us with the discrete problem
\begin{equation}\label{dispro}
\bfgamma = \P_s^\top\Omega \otimes J\, \nabla H\left(\bfe\otimes y_0 + h\I_s\otimes I_{2m} \bfgamma \right).\end{equation}
Once it has been solved, the new approximation (\ref{yh}) is easily seen to be given by
\begin{equation}\label{y1}
y_1 = y_0 + h\gamma_0.
\end{equation}

In \cite{BMR2018} two criteria for {\em a priori} selecting two integer parameters $s_0$ and $s$, $s_0\le s$, are given so that:
\begin{itemize}
\item HBVM$(s_0,s_0)$ solves, up to the round-off error level, the homogeneous linear problem associated with (\ref{odep}),
\begin{equation}\label{odelin}
\dot y = JA y, \qquad y(0)=y_0,
\end{equation}
on the interval $[0,h]$;

\item HBVM$(k,s)$ then solves (\ref{odep}), by using the solution of (\ref{odelin}) to obtain the initial guess for the nonlinear iteration solving (\ref{dispro}). This was indeed paramount, to guarantee its convergence, because of the high-oscillatory nature of the solution;

\item in addition to this, the parameter $k$ defined in (\ref{kappa}) was considered, in order to guarantee (\ref{quadra}).
\end{itemize}
The resulting method was named SHBVM$(k,s,s_0)$ in \cite{BMR2018}, which stands for {\em spectral HBVM with parameters $(k,s,s_0)$}.  

The parameters $s_0$ and $s$ were derived by imposing that the ratio (\ref{ratio}) essentially holds for each frequency component contributing to the solution. It has to be noticed that, according to the analysis in \cite{BMR2018}, the larger the frequencies involved, the larger the parameters $s_0$ and $s$. Consequently, when $\omega$ in (\ref{hop}) is large, $s_0$ and $s$ are large as well. 

The strategy devised in \cite{BMR2018} is finely tuned for {\em highly-oscillatory} systems.  However, when the solution of problem (\ref{odep})--(\ref{hop}) is {\em stiffly-oscillatory}, i.e., only the lowest frequencies contribute to it, whereas the contribution of the highest-frequency components is essentially negligible, the two {\em a priori} criteria defined in \cite{BMR2018} are generally too restrictive, since they would select much larger values for $s$ and $s_0$ than  actually needed. 

In addition, the solution of a  {\em stiffly-oscillatory} problem is often smooth enough that the convergence of the nonlinear iteration for solving (\ref{dispro}) will not require an accurate choice of the initial guess. Consequently,  the solution of the associated homogeneous problem (\ref{odelin}) is no more needed and the role of the parameter $s_0$ becomes quite marginal. As a result, we only need to define an implementation of the HBVM$(k,s)$ method, such that $s$ satisfies
\begin{equation}\label{ratiotol}
\rho_s \equiv \frac{\|\gamma_{s-1}\|}{\max_{i=0,\dots,s-1} \|\gamma_i\|} ~\le~tol\,,
\end{equation}
for a suitably small tolerance $tol\sim u$. This could be in principle done {\em adaptively}, by checking the coefficients $\gamma_j$ at runtime.

In the remaining part of this paper, we shall provide numerical evidence that this can be effectively done for the stiffly-oscillatory Hamiltonian problems deriving from the space semi-discretization of Hamiltonian PDEs, even though we defer to a future paper a thorough analysis for deriving a general criterion.

\section{Space discretization}\label{space}

When solving an Hamiltonian PDE defined in the domain $[a,b]\times[0,T]$, with prescribed initial conditions at $t=0$ and periodic boundary conditions, we shall consider a semi-discretization in space along the Fourier basis:
\begin{eqnarray}\nonumber
c_0(x) &=& \frac{1}{\sqrt{b-a}}, \\  \label{fourier}
c_j(x) &=& \sqrt{\frac{2}{b-a}}\cos\left(2\pi j\frac{x-a}{b-a}\right),\\ \nonumber
 s_j(x) &=& \sqrt{\frac{2}{b-a}}\sin\left(2\pi j\frac{x-a}{b-a}\right),\qquad j=1,2,\dots,
\end{eqnarray}
which is orthonormal, since, for all allowed $i,j$:
$$\int_a^b c_i(x)c_j(x)\dd x = \delta_{ij} = \int_a^b s_i(x)s_j(x)\dd x, \qquad \int_a^b c_i(x)s_j(x)\dd x=0.$$
For simplicity, we confine ourselves to the 1D case, even though the used arguments could in principle be generalized to the case of higher-dimensional space domains. 

The use of the Fourier basis (\ref{fourier}) for the space discretization has been considered in a series of papers \cite{BBFCI2018,BFCI2015,BGS2018,BZL2018} (see also \cite{LIMbook2016,BI2018}), and we collect here some significant examples, i.e., the semilinear wave equation \cite{BFCI2015}, which in first order form reads
\begin{equation}\label{slwe}
u_t = v, \qquad v_t = u_{xx} - f'(u), \qquad (x,t)\in [a,b]\times[0,T],
\end{equation}
and the nonlinear Schr\"odinger equation \cite{BBFCI2018}, which we write in real form as
\begin{equation}\label{nse}
\begin{array}{rl}
u_t &=  -v_{xx}   -f'(u^2+v^2) v, \\[5
pt]
v_t &

=  u_{xx}    +f'(u^2+v^2) u,\qquad (x,t)\in [a,b]\times[0,T].
\end{array}
\end{equation}
Here, $u,v,f$ are real scalar functions, and $f'$ is the derivative of $f$. All the equations are equipped with initial and periodic boundary conditions. The initial conditions, say $u_0(x)$ and $v_0(x)$, as well as $f$, will be assumed to provide a solution which is suitably regular in space (as a periodic function).

For equations (\ref{slwe}) and (\ref{nse}), the solution is expanded in space along the basis (\ref{fourier}), so that, for time dependent coefficients $\alpha_j(t),\beta_j(t),\theta_j(t),\eta_j(t)$, one has
\begin{eqnarray}\label{expu}
u(x,t) &=& \alpha_0(t)c_0(x) + \sum_{j\ge 1} \alpha_j(t)c_j(x) + \beta_j(t)s_j(x), \\ \label{expv}
v(x,t) &=& \theta_0(t)c_0(x) + \sum_{j\ge 1} \theta_j(t)c_j(x) + \eta_j(t)s_j(x).
\end{eqnarray}
The previous expansions can be written in vector form, by introducing the infinite vectors
\begin{equation}\label{wqp}
\bfw(x) = \pmatrix{c} c_0(x) \\ s_1(x)\\ c_1(x)\\ \vdots\endpmatrix, \qquad
\bfq(t) = \pmatrix{c} \alpha_0(t) \\ \beta_1(t)\\ \alpha_1(t)\\ \vdots\endpmatrix, \qquad
\bfp(t) = \pmatrix{c} \theta_0(t) \\ \eta_1(t)\\ \theta_1(t)\\ \vdots\endpmatrix,
\end{equation}
as
\begin{equation}\label{expuv1}
u(x,t) = \bfw(x)^\top\bfq(t), \qquad v(x,t) = \bfw(x)^\top\bfp(t).
\end{equation}
As a result, by introducing the infinite matrix
\begin{equation}\label{DJ}
D = \frac{2\pi}{b-a}\pmatrix{cccc} 0 \\ & 1\cdot J_2\\ &&2\cdot J_2\\ &&&\ddots\endpmatrix, \qquad J_2 = \pmatrix{cc} &1\\ -1\endpmatrix,
\end{equation}
and considering that
$$\int_a^b \bfw(x)\bfw(x)^\top\dd x = I,$$
the identity operator, one verifies that equation (\ref{slwe}) reads
\begin{equation}\label{slwe1}
\dot\bfq = \bfp, \qquad \dot\bfp = -D^\top D\bfq -\int_a^b\bfw(x) f'(\bfw(x)^\top\bfq)\dd x, \qquad t\in[0,T],
\end{equation}
whereas (\ref{nse}) becomes
\begin{eqnarray}\label{nse1}
\dot\bfq &=&  D^\top D\bfp   -\int_a^b \bfw(x) f'\left((\bfw(x)^\top\bfq)^2+(\bfw(x)^\top\bfp)^2\right) \bfw(x)^\top\bfp\,\dd x, \\ \nonumber
\dot\bfp &=&  -D^\top D\bfq   +\int_a^b \bfw(x) f'\left((\bfw(x)^\top\bfq)^2+(\bfw(x)^\top\bfp)^2\right) \bfw(x)^\top\bfq\,\dd x, \qquad t\in[0,T].
\end{eqnarray}
It is quite straightforward to prove the following result.

\begin{theo}\label{Ham1} Problem (\ref{slwe1}) is Hamiltonian with Hamiltonian
\begin{equation}\label{Hslwe1}
H(\bfq,\bfp) = \frac{1}2\left[ \bfp^\top\bfp + \bfq^\top D^\top D\bfq + 2\int_a^b f(\bfw(x)^\top\bfq)\dd x\right].
\end{equation}
Similarly, problem (\ref{nse1}) is Hamiltonian with Hamiltonian
\begin{equation}\label{Hnse1}
H(\bfq,\bfp) = \frac{1}2\left[ \bfp^\top D^\top D \bfp + \bfq^\top D^\top D\bfq -\int_a^b f\left((\bfw(x)^\top\bfq)^2+(\bfw(x)^\top\bfp)^2\right)\dd x\right].
\end{equation}
\end{theo}

\bigskip
We also consider the Korteweg-de Vries equation \cite{BGS2018},
\begin{equation}\label{kdv}
u_t = \nu \, u_{xxx} + \mu \, u u_x, \qquad (x,t)\in [a,b]\times[0,T],
\end{equation}
where $\nu,\mu$ are nonzero real scalars, equipped with initial condition $u_0(x)$ and periodic boundary conditions. The initial condition, as before, is assumed to provide a suitably regular solution in space (as a periodic function). Also the solution of this equation can be expanded along the basis (\ref{fourier}) in the form (\ref{expu}). In such a case, however, it is known that
$$\alpha_0(t)c_0(x) \equiv \int_a^b u_0(x) \dd x \equiv \hat u_0.$$
Consequently, by setting the infinite vectors and matrix
\begin{equation}\label{wyD}
\hat\bfw(x) = \pmatrix{c} c_1(x)\\ s_1(x)\\ c_2(x)\\ s_2(x)\\ \vdots\endpmatrix, \qquad
\bfy(t) = \pmatrix{c} \alpha_1(t) \\ \beta_1(t)\\  \alpha_2(t) \\ \beta_2(t)\\ \vdots\endpmatrix, \qquad
\hat D = \pmatrix{ccc} 1\\ &2 \\ && ~\ddots\endpmatrix,
\end{equation}
the following result can be proved \cite{BGS2018}.
\begin{theo}\label{Ham2} With reference to (\ref{wyD}) and matrix $J_2$ defined in (\ref{DJ}), problem (\ref{kdv}) can be written in Hamiltonian form as
\begin{equation}\label{kdv1}
\dot\bfy = \hat D\otimes J_2\, \nabla H(\bfy)
\end{equation}
with Hamiltonian
\begin{equation}\label{Hkdv1}
H(\bfy) = \frac{1}2\left[  -\nu\left(\bfy^\top \hat D^2\otimes I_2 \bfy\right) +\frac{\mu}3 \int_a^b \left( \hat u_0 + \hat\bfw(x)^\top\bfy\right)^3 \dd x\right].
\end{equation}
\end{theo}

\begin{rem} It can be shown that the Hamiltonian functions (\ref{Hslwe1}), (\ref{Hnse1}), and (\ref{Hkdv1}) are equivalent to the corresponding Hamiltonian functionals defining the corresponding equations \cite{BBFCI2018,BFCI2015,BGS2018}.\end{rem}

For all problems, the Hamiltonian is a constant of motion. For the nonlinear Schr\"odinger equation, there are also the following quadratic invariants \cite{BBFCI2018}, with reference to (\ref{wqp})--(\ref{DJ}):
\begin{equation}\label{M12}
M_1(\bfq,\bfp) = \int_a^b \left[ (\bfw(x)^\top\bfq)^2+(\bfw(x)^\top\bfp)^2\right]\dd x, \qquad M_2(\bfq,\bfp) = 2\left[\bfq^\top D\bfp\right].
\end{equation}

\medskip
In order to derive a numerical method, the expansions (\ref{expu})--(\ref{expv}) need to be truncated at a convenient number $N$ of terms.
In so doing, the vectors (\ref{wqp}) and the matrix (\ref{DJ}) becomes of dimension $2N+1$, whereas the vectors and the matrix in (\ref{wyD}) becomes of dimension $2N$ and $N$, respectively. Upon regularity assumptions on the solution, the truncated expansions (\ref{expu})--(\ref{expv}) converge exponentially to the respective limits, thus providing a {\em spectrally accurate space discretization}. We shall always assume that this will be done, hereafter.

We also mention that the integrals in space, occurring in (\ref{slwe1})--(\ref{Hkdv1}), can be computed (either exactly or approximately within machine precision) by a composite trapezoidal rule based at the points
\begin{equation}\label{xi}
x_i = a+i\frac{b-a}m, \qquad i=0,\dots,m,
\end{equation}
by using a suitably large value of $m$ \cite{BBFCI2018,BFCI2015,BGS2018}.

At last, we mention that the numerical solution of the discrete problem (\ref{dispro}) derived from the application of a HBVM$(k,s)$ to any of the considered problems (\ref{slwe1}), (\ref{nse1}), (\ref{kdv1}) can be made very efficient by using a {\em blended implementation} of the methods \cite{LIMbook2016,BI2018,BIT2011} and considering an approximation of the Jacobian of the right-hand side provided by the linear part only. This latter, in turn, is the same for all time-steps and has a block diagonal structure with diagonal blocks. As a consequence, a very efficient nonlinear iteration can be devised for all the considered PDE problems (we refer to \cite{BBFCI2018,BFCI2015,BGS2018} for full details).

\section{Numerical tests}\label{numtest}

We here compare the following methods:
\begin{itemize}
\item HBVM$(s,s)$ methods, $s=1,2,3$, i.e., the symplectic $s$-stage Gauss methods of order 2,4,6;

\item HBVM$(k,s)$, $s=1,2,3$, and $k$ chosen so that the method is energy-conserving, for the used time-step $\Delta t$;

\item HBVM$(k,s)$, where $s$ is chosen according to (\ref{ratio}) and $k$ according to either (\ref{kappa}) or to gain exact approximation of the integrals (\ref{quadra}). For sake of clarity, we shall refer to such a method as {\em spectral HBVM}.

\end{itemize}
It is worth mentioning that the same code, implemented in Matlab (R\,2017b, running on a 2.8GHz Intel i7 quad-core computer with 16GB of memory), is used for all the above methods. Consequently, the benchmark will be quite homogeneous, both from the software and hardware point of view. All the reported execution times are in seconds.

We shall apply the methods to particular instances of the equations (\ref{slwe}), (\ref{nse}), and (\ref{kdv}) possessing a (known) soliton solution. In all cases, this latter solution is suitably smooth so that the resulting semi-discrete problem is stiffly-oscillatory.

The expansions (\ref{expu})--(\ref{expv}) are truncated at an index $N$  such that the initial conditions are accurately reproduced within a round-off error level, thus providing a spectrally accurate space discretization. Concerning the integration in time provided by a HBVM, choosing $s$ according to the ratio (\ref{ratiotol}) yields a practical spectral accuracy in time. 
The tolerance $tol$ is chosen in order to truncate the expansion (\ref{expy}) when the norm of the last coefficients becomes small and/or ``stagnates'' (meaning that a round-off error level has been reached).
\subsection{Sine-Gordon equation}\label{SGeq}
\begin{figure}[p]
\centerline{\includegraphics[width=9cm,height=5.5cm]{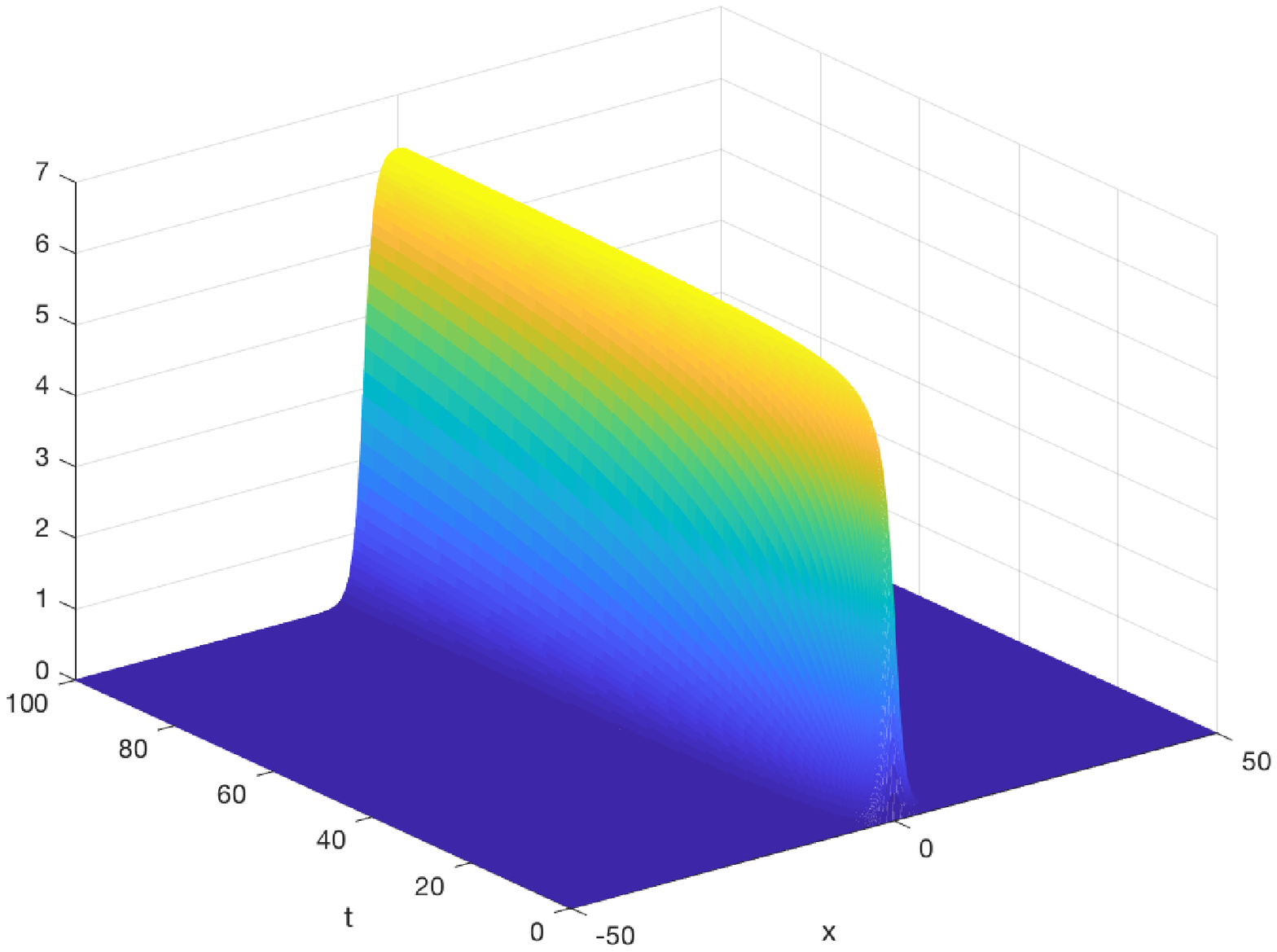}}
\caption{Plot of $u(x,t)$, solution of the Sine-Gordon problem (\ref{doublepole})-(\ref{solSG}).}
\label{SG_fig}
%\end{figure}
%\begin{figure}[t]

\medskip
\centerline{\includegraphics[width=9cm,height=5.5cm]{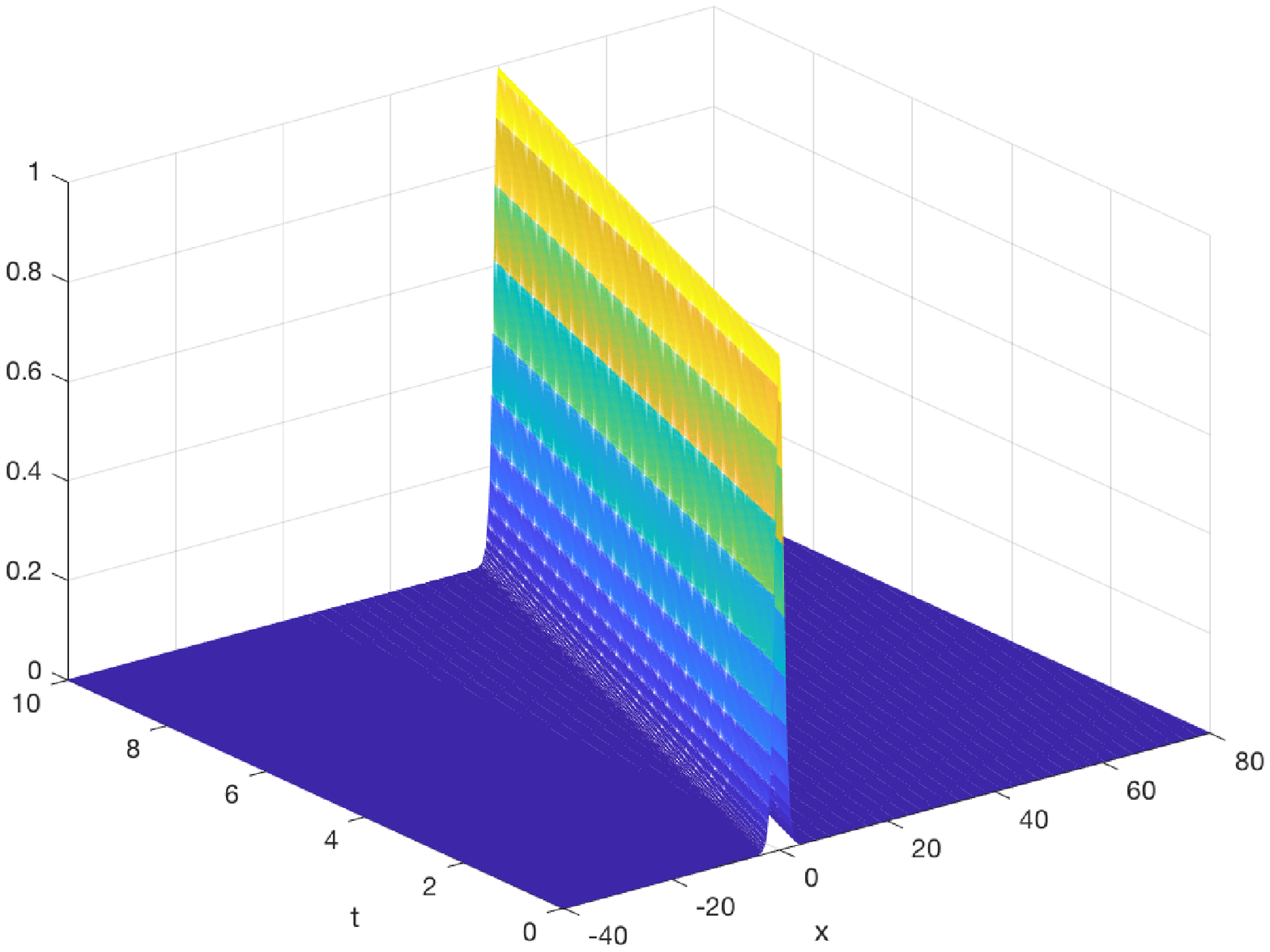}}
\caption{Plot of $u(x,t)^2+v(x,t)^2$, solution of the NLSE problem (\ref{solitone})-(\ref{sol}).}
\label{nlse_fig}
%\end{figure}
%\begin{figure}[t]

\medskip
\centerline{\includegraphics[width=9cm,height=5.5cm]{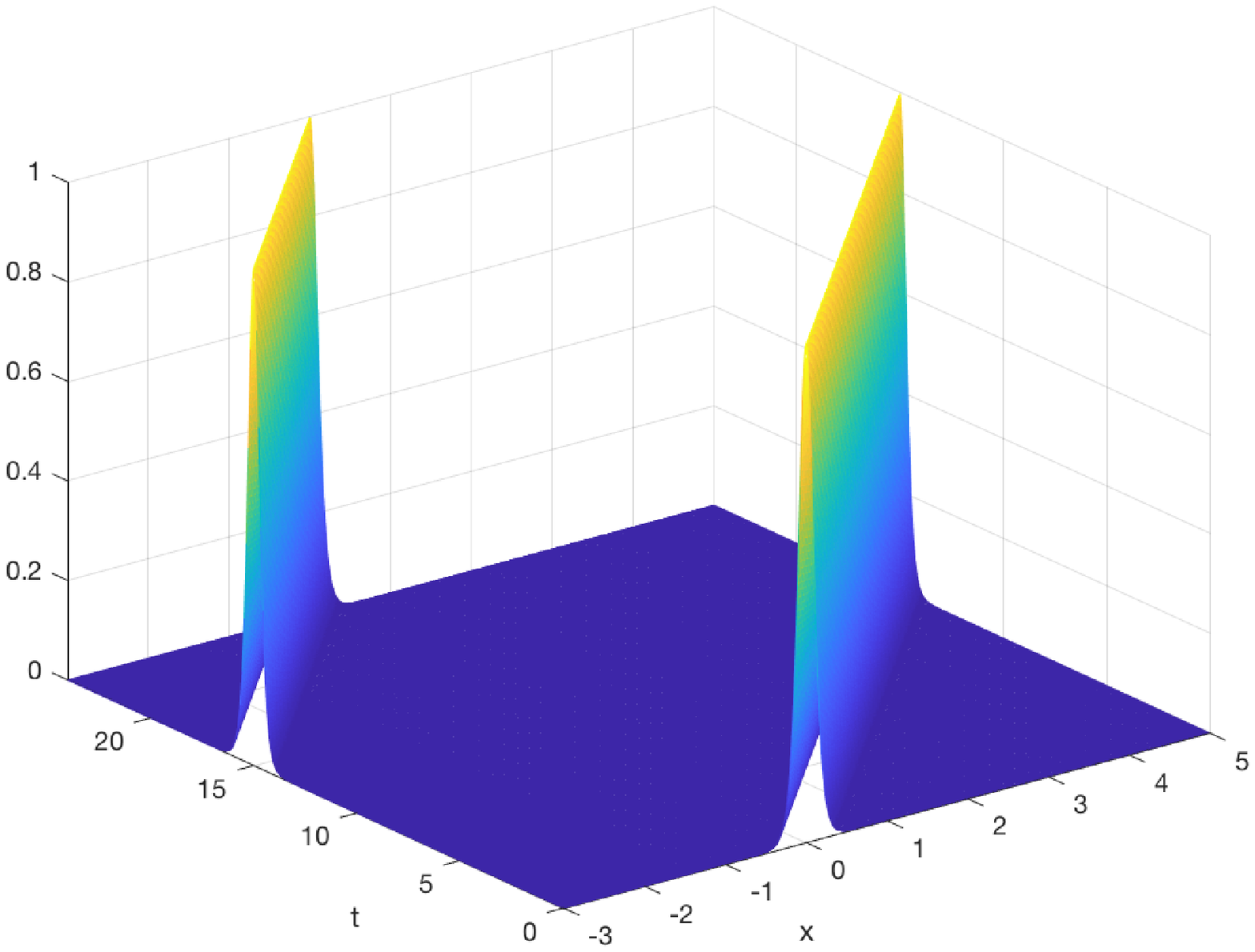}}
\caption{Plot of $u(x,t)$, solution of the KdV problem (\ref{ex1})--(\ref{xiab}).}
\label{kdv_fig}
\end{figure}

In  this example, taken from \cite{BFCI2015},
\begin{equation}\label{doublepole}
u_{tt} = u_{xx} -\sin(u),\qquad     (x,t) \in [-50,50]\times[0,100],
\end{equation}
the initial condition at $t=0$ is obtained from the known solution,
\begin{equation}\label{solSG}
u(x,t) = 4\,\mathrm{atan}(t\,\mathrm{sech}(x)),
\end{equation}
plus periodic boundary conditions. The solution is depicted in Figure~\ref{SG_fig}. In this case, a value $N=300$ has been used for the space discretization, with $m=601$ (see (\ref{xi})) for computing the integrals in space. The corresponding semi-discretization error, measured on the initial condition, is $6.21 \cdot 10^{-14}$. The obtained numerical results are shown in Tables~\ref{SG_gauss}--\ref{SG_shbvm}, where we list the solution error $e_u$ and the Hamiltonian error $e_H$, together with the used time-step and the execution times.
As is expected, from Table~\ref{SG_hbvm} one sees that the energy conserving HBVMs conserve the Hamiltonian.
In Table~\ref{SG_shbvm} we list the results obtained by the spectral HBVM used with parameters $s$ and $k$ computed according to (\ref{ratiotol}) and (\ref{kappa}), with a tolerance $tol$ of the order of $10^{-11}$, able to provide energy conservation and a uniformly small solution error. In the same table, we also report the results obtained by using the original SHBVM$(k^*,s^*,s_0^*)$ method, whose parameters $k^*,s^*, s_0^*$ are computed according to the two criteria given in \cite{BMR2018}. It turns out that both methods have a comparable accuracy but, as was expected, the value $s$ for stiff oscillatory systems is smaller than $s_0^*$ and $s^*$ corresponding to highly oscillatory problems and, consequently, the new method is less time consuming.
 From the obtained results, one infers that the spectral HBVM method here described is the most effective method, with an almost uniform execution time.

\subsection{Nonlinear Sch\"odinger Equation (NLSE)}\label{NLSeq}
In this example, taken from \cite{BBFCI2018},
\begin{eqnarray}\label{solitone}
u_t &=& -v_{xx} - 2(u^2+v^2)v,\\ \nonumber
v_t &=& u_{xx} +2(u^2+v^2)u,\qquad     (x,t) \in [-40,80]\times[0,10],
\end{eqnarray}
the initial condition at $t=0$ is obtained from the known solution,
\begin{equation}\label{sol}
u(x,t) = \mathrm{sech}(x-4t)\cos(2x-3t), \qquad v(x,t) = \mathrm{sech}(x-4t)\sin(2x-3t),
\end{equation}
plus (approximate) periodic boundary conditions. The modulus of the solution $(u^2+v^2)$ is depicted in Figure~\ref{nlse_fig}.
A value $N=300$ has been used for the space discretization, with $m=601$ for computing the integrals in space. The corresponding semi-discretization error is $1.00 \cdot 10^{-14}$. The obtained numerical results are listed in Tables~\ref{nlse_gauss}--\ref{nlse_shbvm}, where we list the solution error $e_{uv}$, the Hamiltonian error $e_H$, the errors on the quadratic invariants (\ref{M12}) ($e_1$ and $e_2$, respectively) together with the used time-step and the execution times. As was expected, the symplectic Gauss methods conserve the two quadratic invariants but not the Hamiltonian (see Table~\ref{nlse_gauss}), whereas the energy conserving HBVMs conserve the Hamiltonian but not the quadratic invariants (see Table~\ref{nlse_hbvm}).
In Table \ref{nlse_shbvm} we give  the parameters $k^*,s^*,s_0^*$ used by the spectral method SHBVM$(k^*, s^*, s_0)$, as defined in \cite{BMR2018}, as well as the parameters $k,s$, obtained by (\ref{kappa}) and (\ref{ratiotol}), used by the spectral HBVM here described, where $tol$ turns out to be of the order of $10^{-14}$, and is able to provide the conservation of all invariants and a uniformly small solution error. It is worth mentioning that, for the larger time-step used, the SHBVM(78,76,40) does not converge at all. Form the above results, one infers that 
 also in this case  the spectral HBVM here described is the most effective method, with a uniformly small execution time.

\subsection{Korteweg-de Vries (KdV) equation}\label{KdVeq}
In this example, taken from \cite{BGS2018},
\begin{eqnarray}\nonumber
u_t(x,t) +  \epsilon u_{xxx}(x,t) + u(x,t)u_x(x,t) &=& 0, \qquad (x,t)\in [-3,5]\times{[0,24]}, \\
\epsilon = 0.0013020833, &&\label{ex1}
\end{eqnarray}
the initial condition at $t=0$ is derived from the known solution of the problem, i.e.,
\begin{equation}\label{uxt1}
u(x,t) = 3c\left[\mathrm{sech}\left( \sqrt{\frac{c}{4\epsilon}}(x-ct)_{[-3,5]} \right)\right]^2, \qquad c=\frac{1}3,
\end{equation}
where, in general,
\begin{equation}\label{xiab}
(\xi)_{[a,b]} := \left\{ \begin{array}{ccl}
\xi, & \mbox{if} & \xi \in[a,b],\\[1mm]
a+\mathrm{rem}(\xi-a,b-a), & \mbox{if} & \xi >b,\\[1mm]
b-\mathrm{rem}(b-\xi,b-a), & \mbox{if} & \xi<a,
\end{array}\right.
\end{equation}
with $\mathrm{rem}$  the remainder in the integer division between the two arguments, plus periodic boundary conditions. As a result, one verifies that the solution (\ref{uxt1}) is periodic in time with period $T=24$. The solution is depicted in Figure~\ref{kdv_fig}.
A value $N=300$ has been used for the space discretization, with $m=901$ for (exactly) computing the integrals in space. The corresponding semi-discretization error is $2.19 \cdot 10^{-14}$. The obtained numerical results are listed in  Tables~\ref{kdv_gauss}--\ref{kdv_shbvm}, where we list the solution error $e_u$ and the Hamiltonian error $e_H$, together with the used time-step and the execution times.
Again, from Table~\ref{kdv_hbvm} one sees that  HBVMs conserve the Hamiltonian function. 
Moreover, in Table \ref{kdv_shbvm} we also list  the parameters $k$ and $s$ provided by (\ref{kappa}) and the ratio (\ref{ratiotol}), where $tol$ turns out to be of the order of $10^{-11}$, and is able to provide energy conservation and a uniformly small solution error.
In this case, the corresponding values of the parameters $k^*,s^*,s_0^*$ for the SHBVM method described in \cite{BMR2018} would be impractically high and, therefore, we do not consider them in the table. Also in this case,  one infers that the spectral HBVM is the most effective method, especially when using the largest time-step.

\section{Conclusions}\label{fine}
In this paper, we have provided numerical evidence that spectral HBVMs, formerly devised for numerically solving highly-oscillatory problems, can be adapted to efficiently handle the stiffly-oscillatory problems deriving from a spectrally accurate space discretization of Hamiltonian PDEs. This is achieved by defining an adaptive strategy to obtain the correct parameters for the method. Numerical tests on the sine-Gordon equation, the nonlinear Schr\"odinger equation, and the Korteweg-de Vries equation duly confirm the effectiveness of the approach, resulting in a spectrally accurate space-time numerical method. It is also worth mentioning that, in principle, this approach could be also used for solving, with spectral accuracy in time, larger classes of problems than that considered here.

\begin{acknowledgements}
The idea of combining spectral accurate discretizations in space and time resulted from interesting discussions of the first author with Volker Mehrmann at the ANODE 2018 Conference.

%The authors wish also to thanks the anonimors referees, for their comments.
\end{acknowledgements}

\begin{table}
\caption{$s$-stage Gauss method (Gauss $s$), for solving the Sine-Gordon problem (\ref{doublepole})-(\ref{solSG}) with time-step $\Delta t = 100/n$.}
\label{SG_gauss}

\smallskip
\centerline{\begin{tabular}{|rrrrrr|}
\hline
\multicolumn{6}{|c|}{Gauss 1}\\
\hline
$n$ & CPU-time & $e_u$ & rate & $e_H$ & rate \\
\hline
 100 &  0.4 & 1.04e\,01 & --- & 1.81e\,00 & --- \\
 200 &  0.5 & 1.14e\,01 & ** & 4.47e-01 & 2.0 \\
 400 &  0.7 & 1.20e\,01 & ** & 1.12e-01 & 2.0 \\
 800 &  1.2 & 1.23e\,01 & ** & 2.79e-02 & 2.0 \\
1600 &  2.2 & 4.27e\,00 & 1.5 & 6.98e-03 & 2.0 \\
3200 &  4.0 & 7.78e-01 & 2.5 & 1.75e-03 & 2.0 \\
6400 &  7.1 & 1.84e-01 & 2.1 & 4.36e-04 & 2.0 \\
12800 & 15.5 & 4.53e-02 & 2.0 & 1.09e-04 & 2.0 \\
25600 & 26.1 & 1.13e-02 & 2.0 & 2.73e-05 & 2.0 \\
\hline
\hline
\multicolumn{6}{|c|}{Gauss 2}\\
\hline
$n$ & CPU-time & $e_u$ & rate & $e_H$ & rate\\
\hline
  100 &  0.9 & 5.55e\,00 & --- & 5.56e-02 & --- \\
  200 &  1.6 & 1.43e\,00 & 2.0 & 8.00e-03 & 2.8 \\
  400 &  2.5 & 6.02e-02 & 4.6 & 3.77e-04 & 4.4 \\
  800 &  3.3 & 3.34e-03 & 4.2 & 2.22e-05 & 4.1 \\
 1600 &  4.9 & 2.02e-04 & 4.0 & 1.37e-06 & 4.0 \\
 3200 &  8.3 & 1.26e-05 & 4.0 & 8.52e-08 & 4.0 \\
 6400 & 14.7 & 7.83e-07 & 4.0 & 5.32e-09 & 4.0 \\
12800 & 28.9 & 4.89e-08 & 4.0 & 3.33e-10 & 4.0 \\
25600 & 57.7 & 3.07e-09 & 4.0 & 2.08e-11 & 4.0 \\
\hline
\hline
\multicolumn{6}{|c|}{Gauss 3}\\
\hline
$n$ & CPU-time & $e_u$ & rate & $e_H$ & rate \\
\hline
  100 &  1.7 & 5.34e-01 & --- & 1.83e-03 & --- \\
  200 &  2.7 & 8.74e-02 & 2.6 & 2.29e-04 & 3.0 \\
  400 &  3.0 & 9.35e-04 & 6.5 & 2.59e-06 & 6.5 \\
  800 &  3.4 & 1.35e-05 & 6.1 & 3.85e-08 & 6.1 \\
 1600 &  5.3 & 2.08e-07 & 6.0 & 5.97e-10 & 6.0 \\
 3200 &  9.8 & 3.23e-09 & 6.0 & 9.27e-12 & 6.0 \\
 6400 & 16.7 & 4.01e-11 & 6.3 & 1.49e-13 & 6.0 \\
12800 & 32.9 & 2.15e-11 & ** & 4.26e-14 & ** \\
\hline
\end{tabular}}
\end{table}

\begin{table}
\caption{Energy-conserving HBVMs, for solving the Sine-Gordon problem (\ref{doublepole})-(\ref{solSG}) with time-step $\Delta t = 100/n$.}
\label{SG_hbvm}

\smallskip
\centerline{\begin{tabular}{|rrrrr|}
\hline
\multicolumn{5}{|c|}{HBVM(5,1)}\\
\hline
$n$ & CPU-time & $e_u$ & rate & $e_H$ \\
\hline
  100 &  1.1 & 1.19e\,01 & --- & 1.21e-07 \\
  200 &  1.4 & 2.99e\,00 & 2.0 & 1.42e-14 \\
  400 &  2.2 & 1.71e-01 & 4.1 & 1.60e-14 \\
  800 &  3.5 & 1.10e-02 & 4.0 & 1.42e-14 \\
 1600 &  6.0 & 7.54e-04 & 3.9 & 1.07e-14 \\
 3200 & 10.6 & 1.76e-04 & 2.1 & 1.60e-14 \\
 6400 & 18.6 & 4.40e-05 & 2.0 & 1.07e-14 \\
12800 & 36.3 & 1.10e-05 & 2.0 & 1.60e-14 \\
25600 & 72.7 & 2.75e-06 & 2.0 & 1.07e-14 \\
\hline
\hline
\multicolumn{5}{|c|}{HBVM(6,2)}\\
\hline
$n$ & CPU-time & $e_u$ & rate & $e_H$  \\
\hline
  100 &  1.2 & 2.81e-01 & --- & 4.85e-08 \\
  200 &  2.3 & 2.31e-03 & 6.9 & 1.42e-14 \\
  400 &  3.6 & 6.55e-05 & 5.1 & 1.07e-14 \\
  800 &  4.9 & 4.24e-06 & 3.9 & 1.42e-14 \\
 1600 &  7.1 & 2.68e-07 & 4.0 & 1.07e-14 \\
 3200 & 12.3 & 1.68e-08 & 4.0 & 1.07e-14 \\
 6400 & 21.7 & 1.05e-09 & 4.0 & 1.07e-14 \\
12800 & 41.2 & 6.56e-11 & 4.0 & 1.07e-14 \\
25600 & 83.8 & 1.27e-11 & ** & 1.60e-14 \\
\hline
\hline
\multicolumn{5}{|c|}{HBVM(6,3)}\\
\hline
$n$ & CPU-time & $e_u$ & rate & $e_H$  \\
\hline
  100 &  2.0 & 2.31e-03 & --- & 1.36e-07 \\
  200 &  3.3 & 1.91e-05 & 6.9 & 7.11e-15 \\
  400 &  3.7 & 3.60e-07 & 5.7 & 7.11e-15 \\
  800 &  4.2 & 5.88e-09 & 5.9 & 1.07e-14 \\
 1600 &  6.7 & 9.28e-11 & 6.0 & 1.07e-14 \\
 3200 & 12.4 & 1.22e-11 & ** & 1.07e-14 \\
 6400 & 21.4 & 1.21e-11 & ** & 1.07e-14 \\
12800 & 41.9 & 1.20e-11 & ** & 1.07e-14 \\
\hline
\end{tabular}}
\end{table}

\begin{table}
\caption{Spectral HBVM$(k,s)$, for solving the Sine-Gordon problem (\ref{doublepole})-(\ref{solSG}) with time-step $\Delta t = 100/n$  (left), along with the obtained result by using the original SHBVM$(k^*,s^*,s_0^*)$ method (right).}
\label{SG_shbvm}

\smallskip
\centerline{\begin{tabular}{|r|rrrrr|rrrrrr|}
\hline
\hline
$n$ & CPU-time & $e_u$ & $e_H$ & $k$ & $s$ &  CPU-time &  $e_u$ &  $e_H$ &  $k^*$ &  $s^*$ &  $s_0^*$  \\
\hline
  100 &  2.4 & 1.23e-11 & 7.11e-15 &  22 &  20 & 3.8 & 5.56e-12 & 8.88e-15 & 38 & 36 & 36 \\
  150 &  1.8 & 1.19e-11 & 7.11e-15 &  20 &  16 & 3.6 & 5.50e-12 & 8.88e-15 & 31 & 29 & 29\\
  200 &  1.9 & 1.26e-11 & 7.11e-15 &  20 &  14 & 3.2 & 5.56e-12 & 8.88e-15 & 28 & 26 & 26\\
\hline
\end{tabular}}
\end{table}

\begin{table}
\caption{$s$-stage Gauss method (Gauss $s$), for solving the NLSE problem (\ref{solitone})-(\ref{sol}) with time-step $\Delta t = 10/n$.}
\label{nlse_gauss}

\smallskip
\centerline{\begin{tabular}{|rrrrrrrr|}
\hline
\multicolumn{8}{|c|}{Gauss 1}\\
\hline
$n$ & CPU-time & $e_{uv}$ & rate & $e_H$ & rate & $e_1$ & $e_2$\\
\hline
 100 &  1.2 & 8.86e-01 & --- & 4.76e-02 & --- & 7.84e-14 & 2.48e-15 \\
 200 &  1.6 & 2.63e-01 & 1.8 & 9.58e-04 & 5.6 & 1.53e-14 & 3.96e-16 \\
 400 &  3.7 & 6.60e-02 & 2.0 & 5.55e-05 & 4.1 & 1.18e-14 & 9.02e-17 \\
 800 &  7.1 & 1.64e-02 & 2.0 & 3.41e-06 & 4.0 & 1.27e-14 & 7.63e-17 \\
1600 & 12.1 & 4.10e-03 & 2.0 & 2.12e-07 & 4.0 & 1.29e-14 & 6.94e-17 \\
3200 & 23.0 & 1.03e-03 & 2.0 & 1.32e-08 & 4.0 & 1.35e-14 & 1.18e-16 \\
\hline
\hline
\multicolumn{8}{|c|}{Gauss 2}\\
\hline
$n$ & CPU-time & $e_{uv}$ & rate & $e_H$ & rate & $e_1$ & $e_2$\\
\hline
 100 &  2.7 & 1.30e-02 & --- & 1.09e-05 & --- & 9.99e-15 & 2.07e-15 \\
 200 &  5.0 & 8.79e-04 & 3.9 & 5.20e-08 & 7.7 & 1.09e-14 & 2.08e-17 \\
 400 &  8.5 & 5.58e-05 & 4.0 & 2.18e-10 & 7.9 & 8.66e-15 & 2.78e-17 \\
 800 & 14.9 & 3.50e-06 & 4.0 & 8.75e-13 & 8.0 & 1.02e-14 & 2.78e-17 \\
1600 & 25.5 & 2.19e-07 & 4.0 & 1.29e-14 & 6.1 & 1.11e-14 & 4.86e-17 \\
3200 & 42.0 & 1.37e-08 & 4.0 & 1.47e-14 & ** & 1.02e-14 & 2.78e-17 \\
\hline
\hline
\multicolumn{8}{|c|}{Gauss 3}\\
\hline
$n$ & CPU-time & $e_{uv}$ & rate & $e_H$ & rate & $e_1$ & $e_2$\\
\hline
 100 &  3.4 & 1.19e-04 & --- & 1.03e-08 & --- & 9.10e-15 & 2.78e-17 \\
 200 &  7.3 & 2.02e-06 & 5.9 & 3.99e-12 & 11.3 & 1.22e-14 & 2.78e-17 \\
 400 & 11.1 & 3.25e-08 & 6.0 & 1.47e-14 & 8.1 & 9.55e-15 & 2.78e-17 \\
 800 & 18.6 & 5.67e-10 & 5.8 & 1.47e-14 & ** & 9.55e-15 & 3.47e-17 \\
1600 & 29.7 & 1.47e-10 & ** & 1.51e-14 & ** & 9.77e-15 & 2.78e-17 \\
\hline
\end{tabular}}
\end{table}

\begin{table}
\caption{Energy-conserving HBVM$(2s,s)$ method, for solving the NLSE problem (\ref{solitone})-(\ref{sol}) with time-step $\Delta t = 10/n$ .}
\label{nlse_hbvm}

\smallskip
\centerline{\begin{tabular}{|rrrrrrrrr|}
\hline
\multicolumn{8}{|c|}{HBVM(2,1)}\\
\hline
$n$ & CPU-time & $e_{uv}$ & rate & $e_H$ & $e_1$ & rate & $e_2$ & rate \\
\hline
 100 &  2.5 & 9.10e-01 & --- & 2.66e-15 & 3.64e-03 & --- & 1.93e-04 & --- \\
 200 &  4.4 & 2.82e-01 & 1.7 & 3.11e-15 & 1.40e-04 & 4.7 & 5.67e-06 & 5.1 \\
 400 &  7.4 & 7.00e-02 & 2.0 & 3.55e-15 & 8.23e-06 & 4.1 & 3.30e-07 & 4.1 \\
 800 & 12.9 & 1.75e-02 & 2.0 & 4.88e-15 & 5.07e-07 & 4.0 & 2.03e-08 & 4.0 \\
1600 & 22.0 & 4.37e-03 & 2.0 & 4.44e-15 & 3.16e-08 & 4.0 & 1.26e-09 & 4.0 \\
3200 & 36.8 & 1.09e-03 & 2.0 & 5.33e-15 & 1.97e-09 & 4.0 & 7.88e-11 & 4.0 \\
\hline
\hline
\multicolumn{8}{|c|}{HBVM(4,2)}\\
\hline
$n$ & CPU-time & $e_{uv}$ & rate & $e_H$ & $e_1$ & rate & $e_2$ & rate \\
\hline
 100 &  2.3 & 1.35e-02 & --- & 3.11e-15 & 1.46e-06 & --- & 5.24e-08 & --- \\
 200 &  5.4 & 8.95e-04 & 3.9 & 4.88e-15 & 6.66e-09 & 7.8 & 2.44e-10 & 7.7 \\
 400 &  9.4 & 5.68e-05 & 4.0 & 3.55e-15 & 2.74e-11 & 7.9 & 1.02e-12 & 7.9 \\
 800 & 15.8 & 3.56e-06 & 4.0 & 4.88e-15 & 1.13e-13 & 7.9 & 4.05e-15 & 8.0 \\
1600 & 27.2 & 2.23e-07 & 4.0 & 4.44e-15 & 1.11e-14 & ** & 1.11e-16 & ** \\
3200 & 44.8 & 1.39e-08 & 4.0 & 4.44e-15 & 1.49e-14 & ** & 1.39e-16 & ** \\
\hline
\hline
\multicolumn{8}{|c|}{HBVM(6,3)}\\
\hline
$n$ & CPU-time & $e_{uv}$ & rate & $e_H$ & $e_1$ & rate & $e_2$ & rate \\
\hline
 100 &  4.2 & 1.20e-04 & --- & 4.44e-15 & 6.88e-10 & --- & 3.48e-11 & --- \\
 200 &  8.1 & 2.03e-06 & 5.9 & 4.00e-15 & 2.52e-13 & 11.4 & 1.30e-14 & 11.4 \\
 400 & 14.0 & 3.26e-08 & 6.0 & 4.00e-15 & 1.09e-14 & 4.5 & 1.11e-16 & 6.9 \\
 800 & 23.0 & 5.65e-10 & 5.9 & 4.88e-15 & 1.22e-14 & ** & 1.04e-16 & ** \\
1600 & 36.9 & 1.47e-10 & ** & 4.44e-15 & 1.18e-14 & ** & 1.32e-16 & ** \\
\hline
\end{tabular}}
\end{table}

\begin{table}
\caption{Spectral HBVM$(k,s)$, for solving the NLSE problem (\ref{solitone})-(\ref{sol}) with time-step $\Delta t = 10/n$  (left), along with the obtained result by using the original SHBVM$(k^*,s^*,s_0^*)$ method (right).}

\label{nlse_shbvm}

\smallskip
\centerline{\begin{tabular}{|rrrrrrrr|rrrrrrrr|}
\hline
\hline
$n$ & CPU & $e_{uv}$ & $e_H$ & $e_1$ & $e_2$ & $k$ & $s$  & CPU & $e_{uv}$ & $e_H$ & $e_1$ & $e_2$ & $k^*$ & $s^*$ & $s_0^*$\\
       & time &                 &            &             &            &        &        & time &   &   &   &   &   &   &  \\
\hline
 100 & 17.0 & 1.06e-10 & 3.55e-15 & 1.09e-14 & 1.04e-16 & 28 & 26 & ** & ** & ** &**& **& 78 & 76 & 40 \\
 150 & 17.8 & 1.06e-10 & 3.11e-15 & 1.18e-14 & 1.11e-16 & 22 & 20 & 44.6 & 1.06e-10 & 3.11e-15 & 8.66e-15 & 9.71e-17 &  61 & 59& 33  \\
 200 & 16.9 & 1.06e-10 & 3.55e-15 & 9.99e-15 & 1.04e-16 & 20 & 16 & 36.6 & 1.06e-10 & 3.55e-15 & 1.22e-14 & 1.25e-16 & 52 & 50 & 29\\
\hline
\end{tabular}}
\end{table}

\begin{table}
\caption{$s$-stage Gauss method (Gauss $s$), for solving the KdV problem (\ref{ex1})-(\ref{xiab}) with time-step $\Delta t = 24/n$.}
\label{kdv_gauss}

\smallskip
\centerline{\begin{tabular}{|rrrrrr|}
\hline
\multicolumn{6}{|c|}{Gauss 1}\\
\hline
$n$ & CPU-time & $e_u$ & rate & $e_H$ & rate \\
\hline
   60 &  3.9 & 1.07e\,00 & --- & 1.12e-02 & --- \\
  120 &  2.3 & 1.01e\,00 & 0.1 & 2.97e-04 & 5.2 \\
  240 &  2.5 & 6.90e-01 & 0.5 & 1.07e-06 & 8.1 \\
  480 &  3.5 & 2.11e-01 & 1.7 & 6.28e-08 & 4.1 \\
  960 &  5.2 & 5.37e-02 & 2.0 & 3.92e-09 & 4.0 \\
 1920 &  8.7 & 1.35e-02 & 2.0 & 2.57e-10 & 3.9 \\
 3840 & 13.3 & 3.37e-03 & 2.0 & 1.54e-11 & 4.1 \\
 7680 & 23.4 & 8.43e-04 & 2.0 & 9.42e-13 & 4.0 \\
15360 & 43.3 & 2.11e-04 & 2.0 & 5.85e-14 & 4.0 \\
\hline
\hline
\multicolumn{6}{|c|}{Gauss 2}\\
\hline
$n$ & CPU-time & $e_u$ & rate & $e_H$ & rate\\
\hline
   60 &  3.8 & 9.12e-01 & --- & 2.51e-03 & --- \\
  120 &  3.9 & 1.39e-01 & 2.7 & 4.80e-04 & 2.4 \\
  240 &  4.9 & 7.61e-03 & 4.2 & 1.90e-05 & 4.7 \\
  480 &  7.7 & 3.00e-04 & 4.7 & 3.49e-09 & 12.4 \\
  960 & 13.0 & 1.86e-05 & 4.0 & 7.39e-12 & 8.9 \\
 1920 & 22.3 & 1.14e-06 & 4.0 & 9.09e-16 & 13.0 \\
 3840 & 37.2 & 7.14e-08 & 4.0 & 5.55e-17 & 4.0 \\
 7680 & 68.9 & 4.45e-09 & 4.0 & 7.63e-17 & ** \\
15360 & 119.5 & 2.78e-10 & 4.0 & 7.63e-17 & ** \\
\hline
\hline
\multicolumn{6}{|c|}{Gauss 3}\\
\hline
$n$ & CPU-time & $e_u$ & rate & $e_H$ & rate \\
\hline
   60 &  3.0 & 1.82e-01 & --- & 7.86e-04 & --- \\
  120 &  4.2 & 2.12e-03 & 6.4 & 2.04e-06 & 8.6 \\
  240 &  6.3 & 5.04e-05 & 5.4 & 3.46e-09 & 9.2 \\
  480 & 10.2 & 1.90e-06 & 4.7 & 1.04e-11 & 8.4 \\
  960 & 18.0 & 5.57e-08 & 5.1 & 1.34e-14 & 9.6 \\
 1920 & 31.2 & 6.17e-10 & 6.5 & 6.25e-17 & 7.7 \\
 3840 & 56.2 & 5.78e-12 & 6.7 & 4.16e-17 & ** \\
 7680 & 94.1 & 8.38e-14 & 6.1 & 1.25e-16 & ** \\
\hline
\end{tabular}}
\end{table}

\begin{table}
\caption{Energy-conserving HBVMs, for solving the KdV problem (\ref{ex1})-(\ref{xiab}) with time-step $\Delta t = 24/n$.}
\label{kdv_hbvm}

\smallskip
\centerline{\begin{tabular}{|rrrrr|}
\hline
\multicolumn{5}{|c|}{HBVM(2,1)}\\
\hline
$n$ & CPU-time & $e_u$ & rate & $e_H$ \\
\hline
   60 & 18.2 & 1.03e\,00 & --- & 1.39e-17 \\
  120 &  6.4 & 9.92e-01 & 0.1 & 1.39e-17 \\
  240 &  7.0 & 5.98e-01 & 0.7 & 1.39e-17 \\
  480 &  8.9 & 1.74e-01 & 1.8 & 1.39e-17 \\
  960 & 13.6 & 4.42e-02 & 2.0 & 1.73e-17 \\
 1920 & 22.3 & 1.11e-02 & 2.0 & 1.73e-17 \\
 3840 & 37.5 & 2.77e-03 & 2.0 & 1.73e-17 \\
 7680 & 65.1 & 6.93e-04 & 2.0 & 2.08e-17 \\
15360 & 108.7 & 1.73e-04 & 2.0 & 2.08e-17 \\
\hline
\hline
\multicolumn{5}{|c|}{HBVM(3,2)}\\
\hline
$n$ & CPU-time & $e_u$ & rate & $e_H$  \\
\hline
   60 &  4.4 & 4.29e-01 & --- & 2.08e-17 \\
  120 &  4.7 & 2.89e-02 & 3.9 & 1.39e-17 \\
  240 &  6.0 & 3.16e-03 & 3.2 & 1.39e-17 \\
  480 &  9.5 & 2.56e-04 & 3.6 & 1.39e-17 \\
  960 & 16.1 & 1.61e-05 & 4.0 & 2.08e-17 \\
 1920 & 27.3 & 9.89e-07 & 4.0 & 2.08e-17 \\
 3840 & 45.8 & 6.20e-08 & 4.0 & 2.08e-17 \\
 7680 & 85.6 & 3.87e-09 & 4.0 & 1.73e-17 \\
15360 & 146.7 & 2.42e-10 & 4.0 & 2.08e-17 \\
\hline
\hline
\multicolumn{5}{|c|}{HBVM(5,3)}\\
\hline
$n$ & CPU-time & $e_u$ & rate & $e_H$  \\
\hline
   60 &  3.1 & 5.40e-02 & --- & 2.08e-17 \\
  120 &  4.1 & 9.18e-04 & 5.9 & 1.39e-17 \\
  240 &  6.2 & 3.00e-05 & 4.9 & 1.73e-17 \\
  480 & 10.6 & 1.01e-06 & 4.9 & 1.39e-17 \\
  960 & 19.1 & 3.06e-08 & 5.1 & 1.39e-17 \\
 1920 & 35.1 & 3.53e-10 & 6.4 & 1.39e-17 \\
 3840 & 60.9 & 3.41e-12 & 6.7 & 2.08e-17 \\
 7680 & 98.6 & 5.14e-14 & 6.1 & 2.08e-17 \\
\hline
\end{tabular}}
\end{table}

\begin{table}
\caption{Spectral HBVM$(k,s)$, for solving the KdV problem (\ref{ex1})-(\ref{xiab}) with time-step $\Delta t = 24/n$.}
\label{kdv_shbvm}

\smallskip
\centerline{\begin{tabular}{|rrrrrr|}
\hline
\hline
$n$ & CPU-time & $e_u$ & $e_H$ & $k$ & $s$ \\ %&\color{blue} $s^{*}$ & \color{blue} $s_0$ \\
\hline
   60 & 12.8 & 3.98e-13 & 1.39e-17 & 20 & 18 \\ %& $*$  & $*$\\
   90 & 17.0 & 9.98e-14 & 1.39e-17 & 20 & 16 \\%  & $*$ & $*$\\
  120 & 20.0 & 4.71e-14 & 1.39e-17 & 20 &14  \\% & $*$ & $*$  \\
\hline
\end{tabular}}
\end{table}

\end{document}